\documentclass[12pt]{amsart}
\usepackage{amscd}
\usepackage{verbatim}
\usepackage[russian,english]{babel}
\usepackage{amssymb, amsmath, amsthm, amscd,ifthen}
\usepackage[dvips]{graphics}
\usepackage[cp866]{inputenc}
\usepackage{graphicx}
\usepackage{epsfig}
\usepackage{hyperref}
\usepackage[all,cmtip]{xy}
\usepackage{tikz-cd}
\usepackage[all]{xy}

\usepackage{amsmath,amssymb,amscd,}
\usepackage[all]{xy}

\usepackage{color}

\usepackage[dvips]{graphics}
\usepackage[cp866]{inputenc}
\usepackage{graphicx}
\usepackage{epsfig}

\theoremstyle{plain}
\newtheorem{theorem}{Theorem}

\theoremstyle{definition}
\newtheorem{definition}{Definition}

\theoremstyle{plain}
\newtoks\thehProclaim
\newtheorem*{Proclaim}{\the\thehProclaim}

\theoremstyle{definition}
\newtoks{\thehRemark}
\newtheorem*{Remark}{\the\thehRemark}

\renewcommand{\leq}{\leqslant}

\begin{document}

	\title{A new proof of Milnor-Wood inequality}

	\author{ Gaiane Panina, Timur Shamazov,  Maksim Turevskii }
	

	

	\begin{abstract} The Milnor-Wood inequality states that if a (topological) oriented circle bundle over an orientable surface of genus $g$ has a
smooth transverse foliation, then the Euler class of the bundle satisfies $$|\mathcal{E}|\leq 2g-2.$$
We give a new proof of the inequality based on  a (previously proven by the authors) local formula which computes $\mathcal{E}$ from the singularities of 
a quasisection.  

We also sketch two other proofs: one based on Poincar\`{e} rotation number theory, and the other of topological nature.
		
	\end{abstract}

	\maketitle

	\section{Introduction}\label{SecIntro}
	
\subsection*{William Thurston and his works on foliations} The doctoral thesis of William Thurston was entitled ''Foliations of Three-Manifolds which are Circle Bundles'', see \cite{CollWorksThurston}. Its main result is very much related to the Milnor-Wood theorem (which had been already known at that time):

\textit{Let $E\rightarrow S$ be an oriented circle bundle over  an oriented surface whose genus is non-zero.  Let $F$ be a transversely  oriented smooth foliation of $E$ without compact leaves homeomorphic to the torus. Then $F$ can be isotoped to a foliation transverse to the fibers.}
 
\bigskip

We refer the reader to  \cite{CollWorksThurston}  for further numerous Thurston's results on foliation and also for \cite{Pap} for a digest.

More than twenty years later,  Thurston wrote, \cite{Thurston1}:
''An interesting phenomenon occurred. Within a couple of years, a dramatic evacuation of the field started to take place. I heard from a number of mathematicians
that they were giving or receiving advice not to go into foliations   they were saying that Thurston was cleaning it out. People told me (not as a complaint, but
as a compliment) that I was killing the field. Graduate students stopped studying
foliations, and fairly soon, I turned to other interests as well.
I do not think that the evacuation occurred because the territory was intellectually exhausted  there were (and still are) many interesting questions that remain
and that are probably approachable. Since those years, there have been interesting
developments carried out by the few people who stayed in the field or who entered
the field, and there have also been important developments in neighboring areas
that I think would have been much accelerated had mathematicians continued to
pursue foliation theory vigorously.
Today, I think there are few mathematicians who understand anything approaching the state of the art of foliations as it lived at that time, although there are some
parts of the theory of foliations, including developments since that time, that are
still thriving.''

\bigskip

\subsection*{What is a circle bundle and its Euler class?}

Formally speaking, an \textit{oriented circle bundle}  over a surface $S_g$ of genus $g$ is a locally trivial fiber bundle
$$ \pi:E\rightarrow S_g$$ whose fibers are oriented circles.

Isomorphism classes of circle bundles over $S_g$ are the same as isomorphism  classes  of oriented $2$-dimensional real linear
vector bundles. They are completely characterized by their Euler classes $\mathcal{E}(E\rightarrow S_g)\in H^2(S_g,\mathbb{Z})$. Since a choice of orientation of $S_g$
identifies $H^2(S_g,\mathbb{Z})$ with $\mathbb{Z}$, we shall speak of the \textit{Euler number} of the bundle.

\bigskip

Informally speaking, a circle bundle with a prescribed Euler number can be constructed as follows.
Take $\pi:S_g\times S^1\rightarrow S_g$, where  $\pi$ is the forgetful projection. This is a trivial bundle, and its Euler number is zero. Next, take a disk $d\subset S_g$; its preimage $\pi^{-1}(d) $ is a solid torus. Cut out $\pi^{-1}(d) $  and patch it back, after applying some power $T^N$ of the Den twist along the fiber. We arrive at another circle bundle whose Euler number equals $N$.

All isomorphisms classes of oriented circle bundles arise in this way, and they can be distinguished by their Euler numbers.

Euler number is an obstruction for existence of continuous section of the bundle.  That is, a continuous section exists iff $\mathcal{E}=0$, which means that  the bundle is trivial.  In other cases there exists a continuous partial section  defined everywhere except a single point. 

\bigskip

\subsection*{What is a transverse foliation?}

A $2$-dimensional foliation of the total space $E$  is a partition of $E$ into $2$-dimensional submanifolds called \textit{leaves},  locally modeled on the decomposition of $\mathbb{R}^3=(x,y,z)$ into all the $2$-planes defined by the equations  $z=const$.

A \textit{transverse foliation} of  a circle bundle $E\xrightarrow[\text{}]{\pi} S_g$  is a smooth foliation of $E$ whose $2$-dimensional {leaves}  are transverse to the fibers of the bundle \cite{Mann}.  Each transverse foliation comes from  a connection form with zero curvature, and vice versa.

\medskip

\textbf{Example. } The trivial bundle $S_g\times S^1\rightarrow S_g$ has a trivial foliation by \textit{horizontal leaves}  $S_g\times \{\alpha\}$.

\medskip

Speaking  at the 17-th International Congress of Mathematics, Thurston \cite{Thurston} said:
''{Given a large supply of some sort of
fabric, what kinds of manifolds can be made from it, in a way that the patterns match up
along the seams? This is a very general question, which has been studied by diverse means
in differential topology and differential geometry. ... A foliation is a manifold made out of
striped fabric, with infinitely thin stripes, having no space between them. The complete
stripes, or leaves, of the foliation are submanifolds; if the leaves have codimension $k$, the
foliation is called a codimension $k$ foliation.}''

\bigskip

Another informal picture of a foliation (represented by a zero-curved connection) is given by D. Sullivan \cite{Sullivan}:
''We assume that the fiber bundle is endowed  by a smooth system of homeomorphisms between the fibers of nearby points $x$ and $x'$ in the base (this means that a { connection} is fixed).  We also assume that for nearby points $x,x',$ and $x''$ a compatibility condition is satisfied  (this means that the curvature of the connection vanishes). This system can be pictured as a transverse foliation.''

\bigskip

\subsection*{Milnor-Wood theorem }

\begin{theorem}\label{MW}
  A (topological) oriented circle bundle $E\xrightarrow[\text{}]{\pi} S_g$ has a
smooth transverse  foliation  if and only if the Euler class of the bundle satisfies $$|\mathcal{E}|\leq 2g-2 .$$
\end{theorem}

\textbf{Remark.} Milnor-Wood inequality appeared initially in a slight disguise in \cite{MilnorIneq} for fiber bundles with the structure group $GL^+(2,\mathbb{R})$ and later on in \cite{Wood}  for all topological bundles, that is, for the structure group $Homeo^+(S^1)$. 
The above formulation  is borrowed from \cite{Mann}, which  gives  an excellent review on the subject and related questions. 

\bigskip
Our main aim is to  give a new proof of the\textit{ Milnor-Wood inequality}, that is, 
 the ''only  if'' part of the theorem.
However let us first shortly comment on the existence of  a foliation  for small values of the Euler number.
 The most important and interesting  is the case $\mathcal{E}=2g-2$. 
Then each foliation of the total space  is \textit{semi-conjugate} (see \cite{Mann}) to the celebrated Anosov foliations that comes from to the geodesic flow.

Here is its construction in short. Let us first construct a foliation of the tangent space to the hyperbolic plane $\mathbb{H}^2$. Given a  point $p$ on the absolute  $\partial \mathbb{H}^2$, take the family of oriented
 geodesics   with the forward limit  at $p$.  The  tangent bundle $T\mathbb{H}^2$ is foliated as follows:
each  leaf consists of the forward unit tangents to a family of
 geodesics with one and the same limit point.  This foliation is invariant under isometries of $\mathbb{H}^2$ and therefore induces a foliation on the tangent space  $TS_g$ whenever $S_g$ is expressed as a quotient of $\mathbb{H}^2$.

\subsection*{Milnor-Wood inequality, a proof coming from Poincar\`{e} rotation number }

\bigskip

 An important observation is that each foliated circle bundle  yields a group homomorphism
  $$\phi:\pi_1(S_g)\rightarrow Homeo^+(S^1).$$
The converse is also true: each such homomorphism yields a foliated bundle.
Here how it goes: cut and unfold $S_g$ to a $4g$-gon  $P$.  Also cut and unfold consistently the total space $E$. After unfolding, we have a trivial circle bundle over
$P$ together with a trivial foliation. To reconstruct the initial fiber bundle, one needs the above homomorphism which gives the   patching rules. 
The number  $\mathcal{E}$  expresses in terms of $\phi$ as the integer Poincar\`{e} rotation number of the product of commutators:
$$\mathcal{E}=rot \prod_{i=1}^g[\phi(a_i),\phi(b_i)] \in \mathbb{Z},$$ 
 where $a_i$ and $b_i$ are the standard generators of $\pi_1(S_g)$.
Now the Milnor-Wood inequality follows by some properties of rotation number.

\bigskip

\subsection*{Milnor-Wood inequality, a topological proof by D. Sullivan  \cite{Sullivan} }

Unfold $S_g$ to a $4g$-gon $P$ in the standard way.  As is known, all the vertices of $P$ patch to a single point $O \in S_g$.

Take a leaf of the foliation. It  gives  a partial section over interior of $P$.  Choose a trivialization of the fiber bundle in a neighborhood of the point $O$, that is, fix an identification of the fibers near $O$ with a single copy of $S^1$. The polygon $P$ has  $4g$ angles. Each of the angles, together with the partial section defines a point on this copy of $S^1$. That is, we may assume that on each of the angles, close to $O$ the section is constant with respect to the chosen trivialization.  Next, we extend the partial section to the edges of $P$ by connecting the corresponding points  by the shortest arcs of $S^1$. This gives a loop $S^1 \rightarrow S^1$ whose degree is the Euler number of the bundle. Since each of the arcs is shorter than $\pi$, the loop is shorter than $4\pi$ and therefore cannot wind more than $2g$ times. Therefore $\mathcal{E}<2g$. It remains to exclude the value $\mathcal{E}=2g-1$. This is explained in the next section.

\section{Milnor-Wood inequality, the new proof}

	The proof  is based  on the \textit{local formula} which  states that Euler number  of a  bundle with the base $S_g$ equals the sum of weights of (some of) singularities  of a quasisection  \cite{PanTurSh}.  The formula is a close relative of the one from \cite{Kaz}.  So we start with a reminder:
	
	\begin{definition}
		A quasisection of $ E \xrightarrow[\text{}]{\pi} S_g$ is an (either bordered or closed) smooth surface  $Q$ and a smooth map $q$
		$$q:Q\rightarrow E$$
		such that $\pi\circ q (Q)=S_g$. 
		The pair $(Q,q)$ will be denoted by $\mathcal{Q}$. We also abbreviate  the composition $\pi\circ q:Q\rightarrow S_g$ as $\pi:{\mathcal{Q}}\rightarrow S_g$.	
	\end{definition}

	We assume  that the maps $q$ and  $\pi\circ q $ are  \textit{generic}, or \textit{stable},  in the sense of Whitney's singularity theory \cite{AVGZ}. In particular, this means that the singularities of $q$ are
	lines of self-intersection, isolated triple points, and Whitney umbrellas only. It is also assumed that away from  Whitney umbrellas, the singularities
	of  $\pi\circ q $  are pleats and folds.

In the present paper we have an even simpler situation: we deal with an embedded disk $\mathcal{Q}$ with    no self-crossings, no folds, no Whitney umbrellas,  no pleats, etc.  The only type of singularities  are transverse self-crossings of  $\pi(\partial \mathcal{Q})$.

\begin{figure}[h]
 \includegraphics[width=0.5\linewidth]{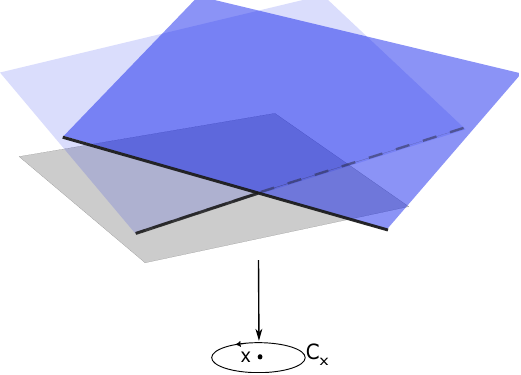}
  \caption {This is the preimage of a neighborhood of a singular vertex. We assume that the fibers are vertical.  Here we have $k=0$,  $n=1$, and $\mathcal{W}_{bb}(x)=1/6$. For its mirror image,  $k=1$,  $n=0$, and $\mathcal{W}_{bb}(x)=-1/6$ }
  \label{Fig1}
\end{figure}

Let $x\in S_g$ be a self-crossing point of  $\pi(\partial \mathcal{Q})$. We call $x$ a \textit{singular vertex} of the quasisection $\mathcal{Q}$.  Locally (in the preimage of a neighborhood of $x$) $\mathcal{Q}$ consists of two bordered sheets, and a non-zero number of regular sheets, see Fig. \ref{Fig1}.

Let $C_x$ be a small circle embracing 	$x$.
	Imagine a point $y$  goes along  the circle $C_x$ in the ccw direction, starting from a place with no bordered sheet in the preimage,
	that is, with the minimal number of points in the preimage $\pi^{-1}(y)\cap\mathcal{ Q}$.
	Let us order the two border lines as follows: the preimage $\pi^{-1}(y)$ meets the first border line   first. The other border line is the second one.
For example,  the right-hand side border line in Figure \ref{Fig1} is the first one.

	A singular vertex
			defines
			two
			numbers, $n$ and $k$:
{ Set
		$n(x)$  be the number of  regular sheets of $\mathcal{Q}$   lying between the first and the second border lines, 
		if one counts from the \textbf{first} border line in the direction of the fiber.}
	Set also $k(x)$ be the number of regular sheets  lying between the border lines, 
	if one counts from the \textbf{second} border line  in the direction of the fiber.

 A singular vertex is assigned a \textit{ weight }by setting
  $$\mathcal{W}_{bb}(x)=\frac{(n-k)}{(n+k)(n+k+1)(n+k+2)}.$$

	Theorem 3 from \cite{PanTurSh} directly implies the following:

	\bigskip

		\textit{ Assume that an embedded quasisection $\mathcal{Q}$ of a circle bundle  has no pleats and no folds. Then the Euler number of the  bundle  equals the sum of weights of 
the singular vertices:
		}
		$$\mathcal{E}= \sum_{x_i} \mathcal{W}_{bb}(x_i).\ \ \ \ \ \ \ (*)$$
		
\bigskip

\begin{figure}[h]
 \includegraphics[width=1.1\linewidth]{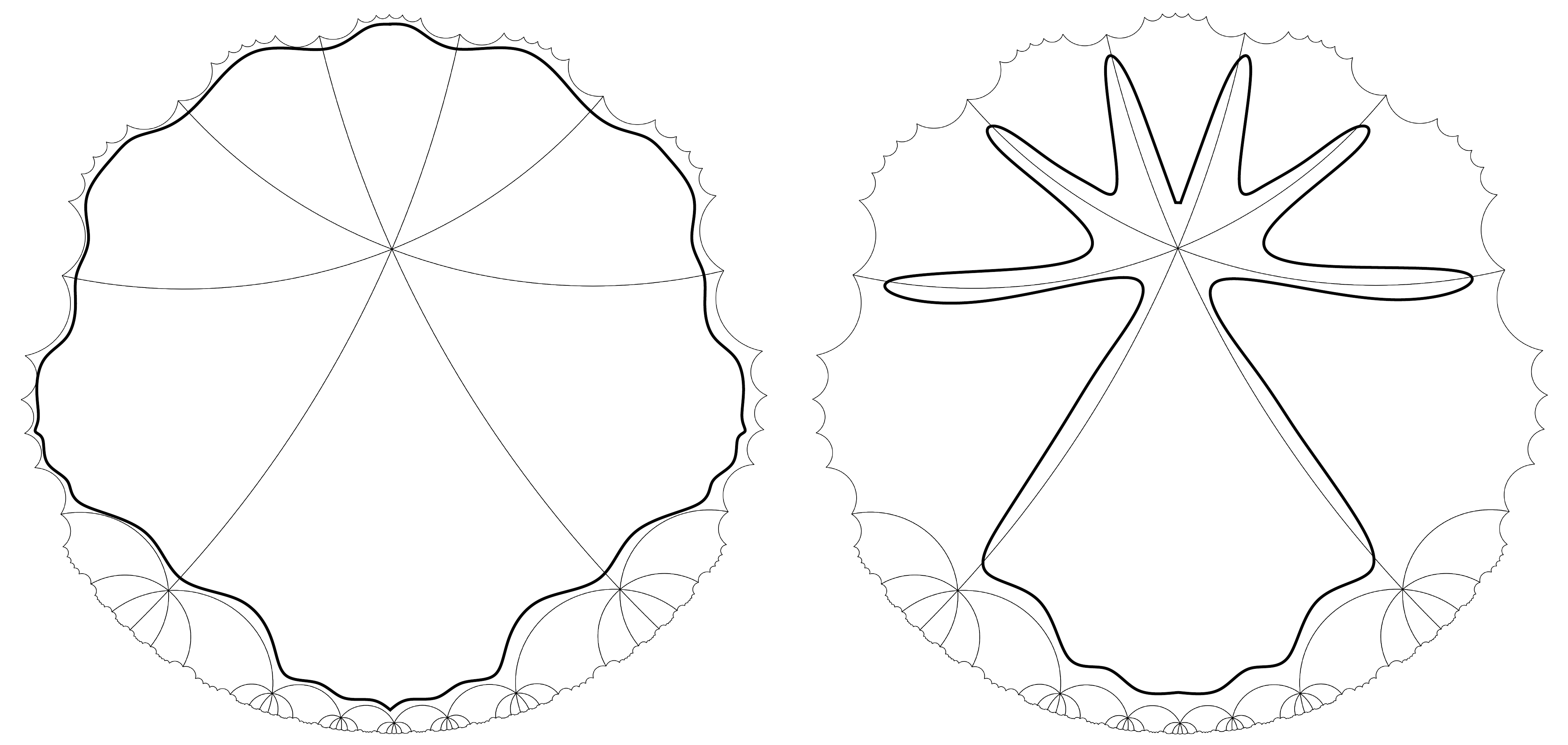}
  \caption {The disks $D_1$ and $D_2$  are bounded by the bold lines. }
  \label{FigDisks}
\end{figure}

	Now we are ready to prove the inequality.
Unfold  $S_g$  to the  regular hyperbolic $4g$-gon  $P$.
The universal cover $U\xrightarrow[\text{}]{pr} S_g$  is tiled by fundamental domains, copies of $P$. 
Consider two discs  $D_1,\ D_2\subset U$, see Fig. \ref{FigDisks}

The disc $D_1$  is slightly smaller than the union of all the fundamental domains that contain some fixed vertex of the tiling.
The disc $D_2$ is contained in $D_1$. It  almost covers  one of the fundamental domains, and has additional ''petals''  so that $pr(D_{1})=S_g$.

\medskip

Let us lift  the disks $D_{1,2}$ to the foliation. This means that we take a map$$\phi_{1,2}:D_{1,2}\rightarrow E $$
 such that  $\phi_{1,2}(D_{1,2})$  lies in a leaf of the foliation, and the below diagram commutes.

\begin{equation}
\label{V_G_SquareDiagram2}
    \xymatrix{
   U\ar[drr]_{pr} &&
D_{1,2} \ar[ll]_{in} \ar[rr]^{\phi_{1,2}} &&  E \ar[lld]^{\pi} &\\
     && S  &\\    }
\end{equation}

We also may assume that $\phi_1(D_1)\subset \phi_2(D_2)$.

The image $\phi_{1,2}(D_{1,2})$ has no self-crossings, but might have overlappings. 

\medskip

\begin{figure}[h]
 \includegraphics[width=0.6\linewidth]{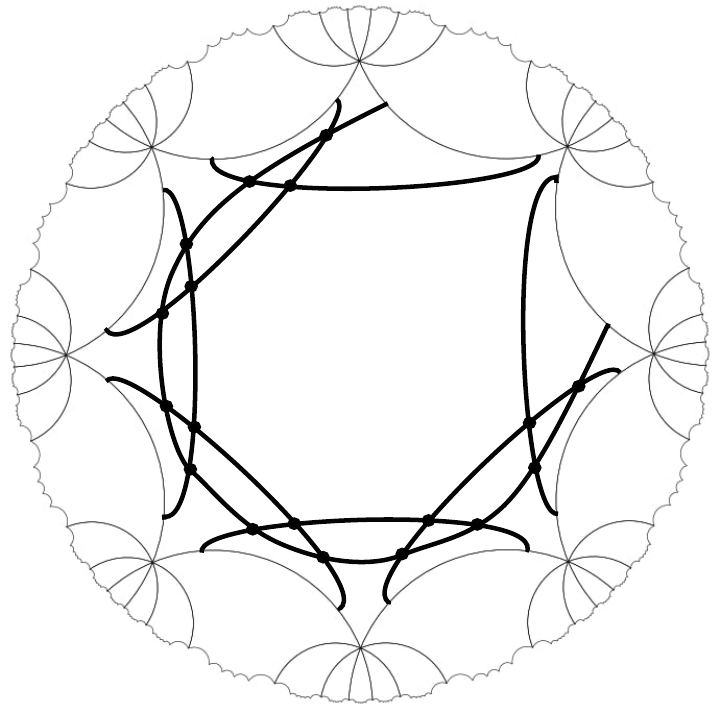}
  \caption{The surface $S_g=S_2$ is unfolded to an octagon. The projection of $\partial \mathcal{Q}$ onto $S_2$ is shown. The bold points indicate the singular vertices. }
  \label{FigSingularities}
\end{figure}

Next, we apply $(*)$ to the quasisection $\mathcal{Q}_1=\phi(D_1)$. While it is possible to list all the singular vertices directly, a simpler approach is to examine the singular vertices  of $\mathcal{Q}_2=\phi(D_2)$,  see Fig.  \ref{FigSingularities} for the case $g=2$.

If  $\phi_2$ is an embedding, there are  $3(4g-2)$ singular vertices,  
each of them contributes $\pm1/6$  since $n+k=1$.
 Their sum  is at most $3(4g-2)\cdot 1/6=2g-1$, so we  have $|\mathcal{E}|\leq 2g-1$.  
If  $\phi_2(D_2)$  has overlaps, some of these singularities disappear, leading to   $|\mathcal{E}|\leq 2g-2$. 

\medskip
It remains to exclude the value $|\mathcal{E}|=2g-1$. We  will do it in two different ways.

\bigskip

\textbf{Excluding the value $|\mathcal{E}|=2g-1$, (1).}

\bigskip

Consider  $\mathcal{Q}_1=\phi_1(D_1)$  and its singular vertices that are the closest to the point $O\in S_g$.
It is sufficient to prove that these   singular vertices  cannot all have the weight $1/6$.

Indeed, the equality is possible only if each of them contributes the maximal possible value of the weight.

Consider an embedded disc $B\subset S_g$  centered at
the point $O$. We assume that $B$ contains the singular vertices that are closest to the point $O$.

\begin{figure}[h]
 \includegraphics[width=0.4\linewidth]{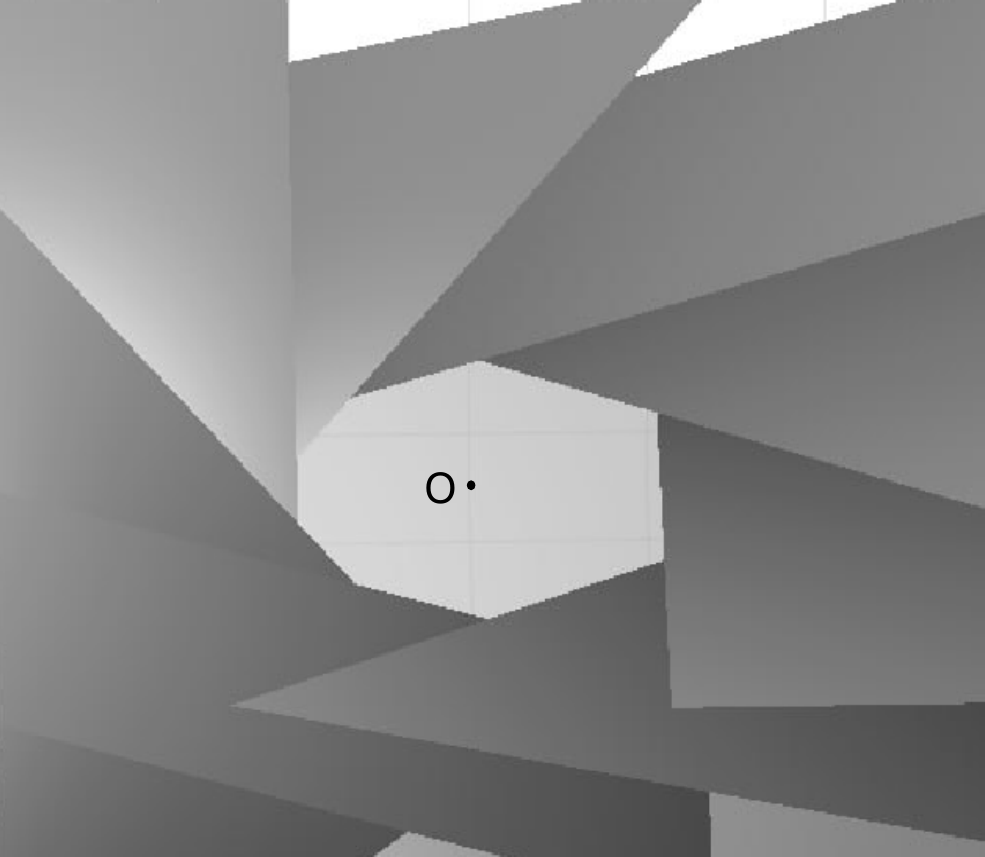}
  \caption {The part of  the quasisection $\mathcal{Q}_1$ lying in $\pi^{-1}(B)$  in assumption that each singular vertex contributes $1/6$. }
  \label{FigEsh}
\end{figure}

Assume that $\mathcal{E}=2g-1$. Then each of these singular vertices  necessarily contributes $1/6$, that is, the local picture is as in Fig. \ref{Fig1}. Thus we  arrive at an ''impossible world of Moris Esher'', see  Fig. \ref{FigEsh}

This means that although the disks $D_2$  $D_2$ embed in $E$ this way, the embeddings do not extend to a transverse foliation (not even to a single leaf).

Let us make the proof more formal.
The restriction of the circle bundle and the restriction  of the foliation to $B$ are trivial,  so we may assume that over $B$, after a suitable coordinatization, we have  horizontal components of $\phi(D_1)$, and therefore,  horizontal $\mathcal{Q}_1$ .

Over $B$, we have $4g$ bordered sheets. Enumerate them counterclockwise as $f_1,...,f_{4g}$.   There is also one regular sheet.
If each of the  singular vertex has the weight $1/6$, the regular sheet cannot lie on the path from $f_1$ to $f_2$ in the direction of the fiber. Analogously, it cannot lie on the path from $f_2$ to $f_3$ in the direction of the fiber, etc. Eventually, there is no room for it. A contradiction.

\bigskip
\textbf{Excluding the value $|\mathcal{E}|=2g-1$, (2).}
\bigskip

This elegant and elementary way of excluding $2g-1$ is borrowed from \cite{Sullivan}.
Take a double cover  $\hat{S}$ of $S_g$; its genus is $2g-1$.  The  foliated bundle $E\rightarrow S_g$ yields  a  foliated bundle over $\hat{S}$ with the Euler number
$\hat{\mathcal{E}}=2\mathcal{E}$. So we have $ |2\mathcal{E}|=|\hat{\mathcal{E}}|\leq 2(2g-1)-1,$ which completes the proof.

	\bigskip
		
		\textbf{Acknowledgement. }  We are indebted to Ilya Alekseev who drew our attention to Milnor-Wood theorem and other beautiful topics around it.
We are thankful to Maxim Kazarian for his interest and a useful remark. We  also thank  Kathryn Mann for pointing out the paper \cite{Sullivan}.

  This work is supported by the Russian Science Foundation (project 25-11-00058)

\end{document}